\documentclass[preprints,article,accept,moreauthors,pdftex]{my-mdpi}

\firstpage{1}
\pubvolume{12}
\issuenum{4}
\articlenumber{391}
\pubyear{2023}
\copyrightyear{2023}
\externaleditor{Academic Editors: Valery Y. Glizer and Feliz Manuel Minhós}
\datereceived{20 January 2023}
\daterevised{31 March 2023}
\dateaccepted{15 April 2023}
\datepublished{18  April 2023}
\hreflink{https://doi.org/10.3390/\\axioms12040391}
\doinum{10.3390/axioms12040391}

\pdfoutput=1 

\Title{Dynamics of a Double-Impulsive Control Model of Integrated Pest Management Using Perturbation Methods and Floquet~Theory} 

\TitleCitation{Dynamics of a Double-Impulsive Control Model of Integrated Pest Management Using Perturbation Methods and Floquet Theory}

\Author{Fahad Al Basir 
$^{1}$\href{https://orcid.org/0000-0003-3744-5524}{\orcidicon}, Jahangir Chowdhury $^{2}$\href{https://orcid.org/0000-0002-9775-9608}{\orcidicon} and Delfim F. M. Torres $^{3,}$*\href{https://orcid.org/0000-0001-8641-2505}{\orcidicon}}

\AuthorNames{Fahad Al Basir, Jahangir Chowdhury and Delfim F. M. Torres}

\AuthorCitation{Al Basir, F.; Chowdhury, J.; Torres, D.F.M.}

\address{$^{1}$ \quad Department of Mathematics, Asansol Girls' College,
Asansol 713304, India\\
$^{2}$ \quad Department of Applied Science, RCC Institute of Information Technology,
Canal South Road,  \mbox{Kolkata 700015, India}\\
$^{3}$ \quad Center for Research and Development in Mathematics and Applications (CIDMA),
Department of Mathematics, University of Aveiro, 3810-193 Aveiro, Portugal}

\corres{Correspondence: delfim@ua.pt; Tel.: +351-234-370-668}

\abstract{We formulate an integrated pest management model to control natural pests
of the crop through the periodic application of biopesticide and chemical pesticides.
In a theoretical analysis of the system pest eradication, a periodic solution
is found and established. All the system variables are proved to be bounded. Our main goal
is then to ensure that pesticides are optimized, in terms of pesticide concentration
and pesticide application frequency, and that the optimum combination of pesticides
is found to provide the most benefit to the crop. By using Floquet theory
and the small amplitude perturbation method, we prove that the pest
eradication periodic solution is locally and globally stable.
The acquired results establish a threshold time limit for the impulsive release
of various controls as well as some valid theoretical conclusions
for effective pest management. Furthermore, after a numerical comparison,
we conclude that integrated pest management is more effective than single
biological or chemical controls. Finally, we illustrate
the analytical results through numerical simulations.}



\begin{document}


\section{Introduction}

In today's farming systems, a variety of approaches are used for pest control.
Maintaining high output while guaranteeing sustainability is crucial for the
entire agriculture sector. Since the beginning of human civilization,
insect and pest control has been one of the most significant difficulties
in the agricultural sector \cite{JTB19,overview19}. Every day, people come up
with fresh ideas for equipment and tactics to use in their fight against pests.
As a result of human efforts to manage pests, our natural ecology
and nature are on the verge of extinction \cite{ref2}.

Chemical controls are less expensive to implement, yet they result in significant
environmental damage \cite{Perring,Bhattacharya,Bhatt}. On the other hand, biological
controls are more costly to implement but have less environmental impact \cite{Fest}.
However, frequently, the use of a single control method is not beneficial to control
pest resistance and preserve environmental quality~\cite{JCV,ML}.
In order to reduce insect populations below economic levels, integrated pest management (IPM),
a safer and more effective method, was developed. IPM is used
for a variety of agronomic crops and is now widely used as an economical and environment-friendly
pest control method in several nations \cite{JCV,ML,Stern}. When the ecological cost of management
is added to the economic price of controls, a combination of chemical and biological controls yields
a superior result when they are used with proper rate and with tolerable intervals. Modeling of
this phenomenon leads naturally to the use of impulsive differential equations \cite{LIU2006,imp17}.

Many researchers have designed mathematical models for pest management
through control strategies, some of which promote chemical agents \cite{Perring,Fest},
some advocate the use of biological agents to impose a total solution of pest and disease
\cite{Wang2010,JTB19,nody22,Stern,overview19,DFMT}, and some researchers use both the
chemical and biopesticides in their mathematical \mbox{models \cite{Bhattacharya,eco_complex,LSC,Chavez}.}
Mathematical-model-based works using impulsive differential equations are also available in the
literature, as already mentioned \cite{Ref_imp1,Ref_imp2,Ref_imp3,Ref_imp4,Ref_imp5,Ref_imp6,Ref_imp7,Ref_imp8}.
Recently, Li, Huang, and Liu proposed a pest management model to simulate
the application of pesticides and build a pesticide function with residual
and delayed effects of pesticides, proposing pest management with pollutant
emission \cite{MR4041509}. Liu et al. constructed a mathematical model
for pest control in which susceptible and infected pests are separated from the pest
population and only susceptible pests are harmful to crops \cite{MR4196616}.
They weighed the two approaches of spraying pesticides and releasing diseased pests
and natural enemies to control vulnerable pests when completing their task.
Alzabut established a mathematical model based on the sense
of biological survey in the field of agriculture, and introduced various
control methods to determine how to protect the crops from destructive pests \cite{Alzabut}.
In \cite{Ref_imp6}, an integrated pest management model using impulsive differential
equations was proposed and analyzed for \textit{Jatropha~curcas} using the release
of infective pests and spraying of chemical pesticides. The existence and stability
of susceptible pest-eradication solutions were analyzed using Floquet theory
and the small amplitude perturbation method. To the best of our knowledge, all
available articles deal with single-impulse differential equation models
where the stability analysis of the periodic pest extinction solution is obtained
by the Floquet theory, the method of small amplitude perturbation,
and the comparison theorem. However, none of the prior research
available in the literature employs the concentration of chemical pesticides
as a system variable as we do here. Moreover, in our study, we spray biological
and chemical pesticides at two different time intervals, simultaneously
varying the time~period.

Pest control models using a single impulse are available (see, e.g.,
\cite{Imp07}), but using double-impulsive controls is rare \cite{eeae_imp,imp09}.
The authors of \cite{eeae_imp} took a predator population along
with biopesticides in an impulsive periodic way for the control
of crop pests. The authors of \cite{imp09} proposed a predator--prey model with disease
in the prey and investigated it for the purpose of integrated pest management. The permanence
of the system and global stability of the susceptible pest-eradication periodic solution
were shown by means of the released amounts of infective prey and predator.
In contrast, here, impulses on both chemical and biopesticides were assumed
in the formulation of the mathematical model for crop pest management.

Our use of the concentration profile of the chemical pesticide as a model
variable is a novel approach. We demonstrated the dynamics using both the
chemical and biological pesticides in the system in an impulsive way, which is,
to the best of our knowledge, a novel concept in crop pest control.
The proposed double-impulsive system was analyzed with proper analytical methods,
namely, using Floquet theory and the perturbation method.

Floquet theory is a powerful mathematical tool for analyzing periodic systems,
and it can be extended to impulsive models with periodic impulses. In impulsive models,
the system's behavior is characterized by a sequence of discrete impulses applied
at regular intervals. These impulses may arise in many practical scenarios,
including electrical circuits, control systems, and biological systems.
Floquet theory provides a robust framework for analyzing and designing
control strategies for impulsive models with periodic impulses. One can analyze
the stability of the impulsive system by examining the eigenvalues of the Floquet matrix \cite{newRef}.
In our analysis, we utilize small amplitude perturbation techniques and Floquet
theory and obtain some valid theoretical results for successful management of pests.
Moreover, we also establish the threshold time limit for the impulsive release of control agents.
Our approach for using the Floquet theory is novel. Additionally, we examined the dynamics
of the system for biological and chemical pesticides used as a sole control measure.

The paper is organized as follows. In Section~\ref{sec2}, we derive the model by using
impulsive differential equations for capturing the IPM system dynamics, taking plant,
pest, virus, and chemical pesticide as model variables. The mathematical analysis of the model
is then discussed in Section~\ref{sec3}, which contains three subsections.
In Section~\ref{subsec3.1}, we determine susceptible pest-eradication periodic solutions
and check the feasibility--boundedness of the system variables discussed in Section~\ref{subsec3.2}.
The local and global stability conditions around the susceptible pest-eradication periodic solutions
are explored in Section~\ref{subsec3.3}. In Section~\ref{sec4}, we exhibit our mathematical results
through numerical simulations. Finally, in Section~\ref{sec5}, we provide a discussion on the three
types of control strategies: spraying chemical pesticide only, impulsively incorporating
of infected pest only, and integrated control with a fixed and a variable impulse
period, to make the final conclusion.


\section{Derivation of the Impulsive Control Model}
\label{sec2}

The following assumptions are taken to formulate the desired model:
the crop plant and susceptible pest populations are denoted by $x$ and $y$,
respectively, and we denote $z$ as the infected pest population.

Due to the finite size of a crop field, which, however, may be large,
we assume logistic growth for the biomass of the crop, with net growth
rate $r$ and carrying capacity $k$. Crops become affected by pests, thereby
causing considerable crop reduction.

Let $\alpha$ be the contact rate between crop and susceptible pest;
let $v(t)$ be the biopesticide (virus); and $s$ be the concentration of chemical pesticide.
A virus infects the susceptible pest at a rate, $\lambda$. The chemical pesticide
kills the susceptible and infected pests at the rates $m_1$ and $m_2$, respectively. Parameters
$c_1$ and $c_2$ are the conversion factors of susceptible and infected pests, respectively,
due to consumption of crop; $d$ and $d+\delta$ are the mortality rates of susceptible and infected pest,
respectively; $\theta$ is the virus replication rate; and $\gamma$ is the lysis rate of the virus.
Finally, we introduce a periodic application of biopesticide and chemical pesticide with different time intervals.

Based on the above assumptions, the desired impulsive system
for integrated pest management is given as
\begin{equation}
\label{mod1}
\begin{cases}
\displaystyle \frac{dx}{dt}
= r x\left(1-\frac{x}{k}\right)-\alpha x y-\phi\alpha x z,
& t\neq (n\tau_1, n\tau_2),\\[0.3cm]
\displaystyle \frac{dy}{dt}
= c_1 \alpha x y-\lambda y v-dy-m_1 s y,
& t\neq (n\tau_1, n\tau_2),\\[0.3cm]
\displaystyle \frac{dz}{dt}
= c_2 \phi \alpha x z+\lambda y v-(d+\delta) z-m_2 s z,
& t\neq (n\tau_1, n\tau_2),\\[0.3cm]
\displaystyle \frac{dv}{dt}
= \theta(d+\delta)z-\gamma v,
& t\neq (n\tau_1, n\tau_2),\\[0.3cm]
\displaystyle \frac{ds}{dt} = -\mu s,
& t\neq (n\tau_1, n\tau_2),\\[0.3cm]
\displaystyle v(t^{+})
=v(t^{-})+ v_i,
&t= n\tau_1,\\[0.3cm]
\displaystyle s(t^{+})
=s(t^{-})+ s_i,
&t= n\tau_2,
\end{cases}
\end{equation}
where $v_i$ and $s_i$ are the strength of biopesticide and chemical pesticide application
in the system at $t= n\tau_1$ and $t= n\tau_2$, respectively; $n = 0, 1, 2, 3,\ldots$,
where $\tau_1$ and $\tau_2$ are the time periods. Here, $v(t^{-})$ and $s(t^{-})$ are the strength
of biopesticide and chemical pesticide before the periodic input, and $v(t^{+})$ and $s(t^{+})$
are the strength of biopesticide and chemical pesticide after the periodic input.

In the impulsive model \eqref{mod1}, we assumed the concentration
of chemical pesticide as a model population, which is realistic and a novel idea.
We use two different impulse intervals for two control agents (biopesticide and chemical pesticide)
that will be analyzed both analytically and numerically. We proceed by analyzing the dynamics
of model \eqref{mod1} by discussing the existence of equilibria with their stability.


\section{Dynamics of the Impulsive Model}
\label{sec3}

In this section, we analyze the boundedness of the solutions of system
(\ref{mod1}), we find out its pest-eradication steady state,
and we analyze the local and global stability. Finally,
we discuss the permanence of the impulsive system.


\subsection{Boundedness of the Model Variables}
\label{subsec3.1}

Let $V(t) = x (t)+ y(t) + z (t)+ v (t)+ s(t)$ and
\begin{eqnarray*}
\frac{dV}{dt} + mV
&\leq&  \big(rx+mx-\frac{r x^2}{k}\big)
-(1-c_1)\alpha x y-(1-c_2)\phi \alpha x z-\{(d-m\}y\\
&&-\{(1-\theta)(d+\delta) -m\} z -(\gamma-m)v-(\mu - m)s.
\end{eqnarray*}

Now, let us define $m = \min\{d, \gamma, \mu, (d+\delta)(1-\theta)\}$.
As $0<c_1, c_2, \theta <1$, then $1-c_1>0$, $1-c_2>0$,
and $1-\theta>0$, so we can write that
\begin{eqnarray}
\frac{dV(t)}{dt} + mV(t) &\leq&M_0,
\end{eqnarray}
where $\displaystyle \frac{k(m+r)^2}{4 r}= M_0$. At $t = n\tau_1$, we have
\begin{eqnarray}
V(n\tau_1^+)\leq V(n\tau_1) + v_i.
\end{eqnarray}
By the comparison theorem, for $t\geq0$ we have
\begin{eqnarray}
V(t)&\leq& V(0)e^{-mt} + \frac{M_0(1-e^{-mt})}{m}+v_i
\frac{e^{-m(t-\tau_1)}}{1-e^{m\tau_1}}+v_i\frac{e^{m\tau_1}}{e^{m\tau_1}-1}\nonumber\\
&\rightarrow& \frac{M_0}{m}+v_i\frac{e^{m\tau_1}}{e^{m\tau_1}-1}~\mbox{as}~t\rightarrow \infty.
\end{eqnarray}
When $t = n\tau_2$,
\begin{eqnarray}
V(n\tau_2^+)\leq V(n\tau_2) + v_i
\end{eqnarray}
and from the comparison theorem it follows that
\begin{eqnarray}
V(t)&\leq& V(0)e^{-mt} + \frac{M_0(1-e^{-mt})}{m}
+v_i\frac{e^{-m(t-\tau_2)}}{1-e^{m\tau_2}}+v_i\frac{e^{m\tau_2}}{e^{m\tau_2}-1}\nonumber\\
&\rightarrow& \frac{M_0}{m}+v_i\frac{e^{m\tau_2}}{e^{m\tau_2}-1}~\mbox{as}~t\rightarrow \infty.
\end{eqnarray}
Thus, $V(t)$ is uniformly bounded
and there exists a positive constant $M>0$ such that
$x(t)\leq M$, $y(t)\leq M$, $z(t)\leq M$, $v(t)\leq M$, and $s(t)\leq M$ for all $t$.

From the above discussion, we have the following theorem.

\begin{Theorem}
\label{thm1}
For the impulsive system (\ref{mod1}), there exists a positive constant $M$
such that $x(t)\leq M$, $y(t)\leq M$, $z(t)\leq M$, $v(t)\leq M$,
and $s(t)\leq M$ for all $t$.
\end{Theorem}

For non-negative solutions, the following lemma follows from \cite{Wang2010}.

\begin{Lemma}
Let $X(t)$ be a solution of the impulsive system (\ref{mod1})
with $X(0^+)\geq 0$. Then $X(t)\geq 0$ for all $t > 0$.
\end{Lemma}


\subsection{Existence of the Pest-Free Periodic Orbit}
\label{subsec3.2}

Since both pests are assumed to be harmful for crops, we discuss
stability at infected and susceptible pest-eradication solutions
of the system when $y=0$ and $z=0$, $t\neq (n\tau_1, n\tau_2)$,
and the linear forms of the fourth and fifth equation of (\ref{mod1}) are
\begin{eqnarray}
\label{mod3}
\frac {dv}{dt} = -\gamma v~~\mbox{and~~}\frac {ds}{dt} =  -\mu s,
\end{eqnarray}
respectively. For an impulse control, we must have
\begin{equation}
\label{mod4}
\begin{split}
v(t^{+}) &= v(t^{-})+v_i,~\mbox{ for }~t= \tau_1,\\
s(t^{+}) &= s(t^{-})+s_i,~\mbox{ for }~t= \tau_2.
\end{split}
\end{equation}
From (\ref{mod3}) and (\ref{mod4}) it is clear that $v$ and $s$
are independent of each other. Thus, the solution of Equation
(\ref{mod3}) can be given as follows:
\begin{eqnarray}
\label{mod5}
v(t) &=&  \{v(0^{+})-v^*(0^{+})\}e^{-\gamma t}+v^*(t),~for~ t\in(\tau_1, (n+1)\tau_1], \nonumber\\
s(t) &=&  \{s(0^{+})-s^*(0^{+})\}e^{-\mu t}+s^*(t), ~for~\in(\tau_2, (n+1)\tau_2],
\end{eqnarray}
where $v^*(t)$ and $s^*(t)$, the positive periodic solution of (\ref{mod3}), are given by
\begin{eqnarray}
\label{mod6}
v^*(t) &=& \frac{v_i e^{-\gamma (t-n\tau_1)}}{1-e^{-\gamma \tau_1}},
~~~~ s^*(t) = \frac{s_i e^{-\mu (t-n\tau_2)}}{1-e^{-\mu \tau_2}},
\end{eqnarray}
with initial values
\begin{equation}
\label{mod7}
v^*(0^+) = \frac{v_i }{1-e^{-\gamma \tau_1}},
\quad s^*(0^+) = \frac{s_i}{1-e^{-\mu \tau_2}}.
\end{equation}
If $y(t)=0$ and $z(t)=0$, then the first equation of (\ref{mod1}) is
\begin{eqnarray}
\label{mod8}
\frac {dx}{dt} &=& r x\left(1-\frac{x}{k}\right),
\end{eqnarray}
which is a logistic equation, and its solution is
\begin{equation}
\label{mod9}
x(t) = \frac{k x(0)}{x(0)+(k-x(0))e^{r t}}~
\text{ for } ~t\neq (n\tau_1, n\tau_2).
\end{equation}
Clearly, (\ref{mod9}) has two equilibria, such as $x=0$ and $x=k$. Therefore, (\ref{mod1})
has two pest-eradication solutions, $(0,0,0,v^*,s^*)$ and $(k,0,0,v^*,s^*)$. Obviously,
at $x=0$ the system (\ref{mod9}) is impossible from the perspective of ecology. For this reason,
in the following subsection we study the stability for the system (\ref{mod1}) at $E=(k,0,0,v^*,s^*)$.


\subsection{Stability of the Pest-Free Periodic Solution}
\label{subsec3.3}

We establish the following theorem for the stability of the pest-free periodic orbit.

\begin{Theorem}
System (\ref{mod1}) is both locally and globally stable
around the pest-free periodic solution $E=(k,0,0,v^*,s^*)$ for the following:
\begin{itemize}
\item[(i)] Application of biopesticide and chemical pesticide
with same time interval $t = n\tau$, provided that
\begin{eqnarray}\label{cond1}
&& c_1 \alpha -d-\frac{\lambda v_i e^{-\gamma (t-n\tau)}}{1-e^{-\gamma \tau}}
- \frac{m_1 s_i e^{-\mu (t-n\tau)}}{1-e^{-\mu \tau}}<0, \nonumber\\
&& c_2 \phi k \alpha -(d+\delta) - \frac{m_2 s_i e^{-\mu (t-n\tau)}}{1-e^{-\mu \tau}}<0;
\end{eqnarray}

\item[(ii)] Application of biopesticide with time interval $t = n\tau_1$ and
chemical pesticide with time interval $t = n\tau_2$, i.e., for different
time intervals, where $\tau_1\neq \tau_2$, provided that
\begin{equation}
\label{cond3}
\begin{gathered}
c_1 \alpha -d-\frac{\lambda v_i
e^{-\gamma (t-n\tau_1)}}{1-e^{-\gamma \tau_1}} <0,\\
c_2 \phi k \alpha -(d+\delta) <0,\\
c_1 \alpha -d- \frac{m_1 s_i e^{-\mu (t-n\tau_2)}}{1-e^{-\mu \tau_2}}<0.
\end{gathered}
\end{equation}
\end{itemize}
\end{Theorem}

\begin{proof}
We need to prove the stability of the system in two cases:

(i) Application of chemical pesticide and biopesticide with same time interval;

(ii) Application of biopesticide and chemical pesticide with different time intervals.

(i) In this case, let $t= nt_1=nt_2=n\tau$. 
We discuss the stability
of the system through the small amplitude perturbation method
at the periodic solution $(k,0,0,v^*,s^*)$. Let
\begin{eqnarray}
\label{mod10}
&&x(t) = k +\epsilon_1(t), \quad y(t)=\epsilon_2(t), \quad z(t)=\epsilon_3(t),\nonumber\\
&&v(t)=v^*(t)+\epsilon_4(t), \quad s(t)=s^*(t)+\epsilon_5(t).
\end{eqnarray}
Here, $\epsilon_1, \epsilon_2, \epsilon_3, \epsilon_4$, \mbox{and}
$\epsilon_5$ denote small amplitude perturbations. Thus, the corresponding
system of (\ref{mod1}) at $(k,0,0,v^*,s^*)$ is given by
\begin{eqnarray}
\label{mod2}
\frac {d\epsilon_1}{dt}
&=& r \{k +\epsilon_1(t)\}\left(1-\frac{k +\epsilon_1(t)}{k}\right)
-\alpha \{k +\epsilon_1(t)\} \epsilon_2(t)\nonumber\\
&&-\phi \alpha \{k +\epsilon_1(t)\} \epsilon_3(t),
~~~~~~~~~~~~~~~~~~~~~~~~~~~~~~~~~t\neq n\tau,\nonumber\\
\frac {d\epsilon_2}{dt} &=& c_1 \alpha \{k +\epsilon_1(t)\}
\epsilon_2(t)-\lambda \epsilon_2(t) \{v^*(t)+\epsilon_4(t)\}\nonumber\\
&&-d\epsilon_2(t)-m_1\{s^*(t)+\epsilon_5(t)\}\epsilon_2(t),
~~~~~~~~~~~~~~~~~t\neq n\tau ,\nonumber\\
\frac {d\epsilon_3}{dt} &=& c_2 \phi \alpha \{k +\epsilon_1(t)\}
\epsilon_3(t)+\lambda \epsilon_2(t) \{v^*(t)+\epsilon_4(t)\}
-(d+\delta) \epsilon_3(t)\nonumber\\
&&-m_2\{s^*(t)+\epsilon_5(t)\}\epsilon_3(t),
~~~~~~~~~~~~~~~~~~~~~~~~~~~~~~~t\neq n\tau,\nonumber\\
\frac {d\epsilon_4}{dt} &=& \theta (d+\delta) \epsilon_3(t)
-\gamma \{v^*(t)+\epsilon_4(t)\},
~~~~~~~~~~~~~~~~~~~~~t\neq n\tau,\nonumber\\
\frac {d\epsilon_5}{dt} &=& -\mu \{s^*(t)+\epsilon_5(t)\},
~~~~~~~~~~~~~~~~~~~~~~~~~~~~~~~~~~~~~~~t\neq n\tau,\nonumber\\
\triangle \{v^*(t)+\epsilon_4(t)\}&=& v_i,
~~~~~~~~~~~~~~~~~~t= n\tau,\nonumber\\
\triangle \{s^*(t)+\epsilon_5(t)\}&=& s_i,
~~~~~~~~~~~~~~~~~~t= n\tau.
\end{eqnarray}
Now, the linear system corresponding to the system \eqref{mod2} is given as
\begin{eqnarray}
\label{mod13}
\frac{d\epsilon_1}{dt}
&=&-r\epsilon_1(t)-\alpha k \epsilon_2(t)-\phi \alpha k\epsilon_3(t),
~~~~~~~~~~~~~~~~~~~~~~~~~~~~~~~~~~~~~~~~~t\neq n\tau,\nonumber\\
\frac{d\epsilon_2}{dt}
&=&c_1 \alpha k \epsilon_2(t)-\lambda \epsilon_2(t)v^*(t)
-d \epsilon_2(t)-m_1 s^*(t)\epsilon_2(t),
~~~~~~~~~~~~~~~~~~ t\neq n\tau,\nonumber\\
\frac{d\epsilon_3}{dt}
&=&c_2 \phi k \alpha \epsilon_3(t)+\lambda
\epsilon_2(t)v^*(t)-(d+\delta) \epsilon_3(t)-m_2 s^*(t)\epsilon_3(t),
~~ ~~~~~ t\neq n\tau,\nonumber\\
\frac{d\epsilon_4}{dt}
&=&\theta (d+\delta) \epsilon_3(t)-\gamma \epsilon_4(t),
~~~~~~~~~~~~~~~~~~~~~~~~~~~~~~~~~~~~~~~~~~~~~~~~~t\neq n\tau,\nonumber\\
\frac{d\epsilon_5}{dt}
&=&-\mu\epsilon_5(t),
~~~~~~~~~~~~~~~~~~~~~~~~~~~~~~~~~~~~~~~~~~~~~~~~~~~~~~~~~~~~~~~~~~
t\neq n\tau,\nonumber\\
\triangle \epsilon_4(t)
&=& v_i, \quad t= n\tau,\nonumber\\
\triangle \epsilon_5(t)
&=& s_i, \quad t= n\tau.
\end{eqnarray}
The fundamental matrix $M(t)$ of (\ref{mod13}) is obtained as
$$
\frac{dM(t)}{dt}=\frac{1}{m}\left[\begin{array}{ccccc}
-r & -\alpha k & -\phi \alpha k& 0 & 0\\
0 & c_1 \alpha k-\lambda v^*(t)-d -m_1 s^*(t) & 0 & 0 & 0\\
0 & \lambda v^*(t) & m_{33}& 0 & 0\\
0 & 0 & \theta(d+\delta) & -\gamma & 0\\
0 & 0 & 0 & 0 & -\mu\\
\end{array}\right]
$$
with initial condition $M(t)=I_5$ (the identity matrix) and
$m_{33}=c_2 \phi k \alpha -(d+\delta) -m_2 s^*(t)$.
Now, the fundamental solution matrix is given by
$$
M(t)=\frac{1}{m}\left[
\begin{array}{ccccc}
\exp{(-r t)} & M_1(t) & M_2(t) &  0 & 0\\
0 & M_3(t) & 0 & 0 & 0\\
0 & M_4(t) &  M_5 & 0 & 0\\
0 & 0 & M_5(t) & \exp{(-\gamma t)} & 0\\
0 & 0 & 0 & 0 & \exp{(-\mu t)}\\
\end{array}\right].
$$

Here, $M_5=\exp\displaystyle \int_0^t\{c_2 \phi k \alpha
-(d+\delta) -m_2 s^*(t)\}dt$,
\begin{eqnarray}
M_3(t)&=&\exp\int_0^\tau\{c_1 \alpha k
-\lambda v^*(t)-d -m_1 s^*(t)\}dt, \nonumber
\end{eqnarray}
where the other $M_i(t)$s are not required for our further analysis.
According to Floquet theory \cite{newRef},
the periodic solution $E(k,0,0,v^*v,s^*)$ is
asymptotically stable if the absolute values
of the eigenvalues of $M(\tau)$ are less than one.

The eigenvalues of $M(t)$ are
\begin{eqnarray}
&&\lambda_1=\exp\{-r \tau\},
\lambda_2=\exp\int_0^\tau\{c_1 \alpha k-\lambda v^*(t)-d -m_1 s^*(t)\}dt, \nonumber\\
&&\lambda_3=\exp\int_0^\tau\{c_2 \phi k \alpha -(d+\delta) -m_2 s^*(t)\}dt,
\lambda_4=\exp\{-\gamma \tau\}, \nonumber\\
&&\lambda_5=\exp\{-\mu \tau\}.\nonumber
\end{eqnarray}
Clearly, $0<\lambda_1<1$, $0<\lambda_4<1$ and $0<\lambda_5<1$. Thus, when both pesticides are
applied with the same time interval, then the system is locally stable around the periodic solution
$E=(k,0,0,v^*v,s^*)$ if  $0<\lambda_2<1$ and $0<\lambda_3<1$. From this, we obtain
\begin{eqnarray}\label{LS0}
&& \exp\int_0^\tau\{c_1 \alpha k
-\lambda v^*(t)-d -m_1 s^*(t)\}dt<1, \nonumber\\
&& \exp\int_0^\tau\{c_2 \phi k \alpha -(d+\delta) -m_2 s^*(t)\}dt<1.
\end{eqnarray}
From Equation (\ref{LS0}), we can choose $\delta_1>0$ such that
\begin{eqnarray}
&& \eta_1=\exp\int_{n\tau}^{(n+1)\tau}\{c_1 \alpha
-\lambda (v^*(t)-\delta_1)-d -m_1 (s^*(t)-\delta_1)\}dt<1, \nonumber\\
&& \eta_2=\exp\int_{n\tau}^{(n+1)\tau}\{c_2 \phi k \alpha
-(d+\delta) -m_2 (s^*(t-\delta_1))\}dt<1.\nonumber
\end{eqnarray}
Since all state variables are positive,
\begin{eqnarray*}
\frac {dv}{dt} \geq-\gamma v;
\quad \frac {ds}{dt} \geq -\mu s.
\end{eqnarray*}
Thus, according to the comparison theorem and Equations (\ref{mod6}) and (\ref{mod7}),
for $\delta_1>0$ there exists $t_0>0$ such that $v(t)\geq v^* - \delta_1$, $s(t)
\geq s^* - \delta_1$ for all $t>t_0$.

From the second equation of system (\ref{mod1}), it can be written that
\begin{eqnarray}
\label{mod15}
\dot{y}(t) &\leq& y(t)\{c_1 \alpha -\lambda (v^*(t)-\delta_1)
-d -m_1 (s^*(t)-\delta_1)\},
\quad t\neq n\tau,\nonumber\\
y(t^+)&=& y(t),
\quad t= n\tau.
\end{eqnarray}
Integrating (\ref{mod15}) into $[n\tau, (n+1)\tau]$, it can be shown that
\begin{eqnarray}
\label{mod16}
y\{(n+1)\tau\} &\leq& y(n\tau)\exp\int_{n\tau}^{(n+1)\tau}
\{c_1 \alpha -\lambda (v^*(t)-\delta_1)-d -m_1 (s^*(t)-\delta_1)\}dt\nonumber\\
&=& y(n\tau)\eta_1.
\end{eqnarray}
Similarly,
\begin{eqnarray}
\label{mod17}
y\{n\tau\} &\leq& y\{(n-1)\tau\}\eta_1.
\end{eqnarray}
Hence, from (\ref{mod16}) and (\ref{mod17}),
\begin{eqnarray*}
y\{(n+1)\tau\} &\leq& y\{(n-1)\tau\}\eta_1^2.
\end{eqnarray*}

Proceeding in this way, we obtain
\begin{eqnarray}
\label{mod18}
y\{(n+1)\tau\} &\leq& y(\tau)\eta_1^n.
\end{eqnarray}
Since $\eta_1<1$, one has $\eta_1^n\rightarrow 0$
whenever $n\rightarrow \infty$. Hence,
$y\{(n+1)\tau\}\rightarrow 0$ as $n\rightarrow \infty$.
Now we take $n\tau<t~\leq~(n+1)\tau$. Then, clearly,
$0<y(t)~\leq~y(n\tau)~\exp~{( n~\tau)}$. Thus,
$y(t)\rightarrow 0$ as $t\rightarrow \infty$.
For
\begin{eqnarray}
\eta_2=\exp\int_{n\tau}^{(n+1)\tau}\{c_2 \phi k \alpha
-(d+\delta) -m_2 (s^*(t-\delta_1))\}dt<1,\nonumber
\end{eqnarray}
we can similarly prove that $z(t)\rightarrow 0$ as $t\rightarrow \infty$.

We now prove that $v(t)\rightarrow v^*(t)$ as $t\rightarrow \infty$.
Since $z(t)\rightarrow 0$  as $t\rightarrow \infty$, then
for some $0<\delta_2<\frac{\gamma}{\theta(d+\delta)}$ there exists
$t_1>0$ such that $0<z(t)<\delta_2$ for all $t>t_1$. Thus, for $t>t_1$
and from the fourth equation of system (\ref{mod1}), we can write that
\begin{eqnarray}
\label{mod19}
\theta(d+\delta)\delta_2 -\gamma v(t)\geq\dot{v}(t)
\geq -\theta(d+\delta)\delta_2 -\gamma v(t).
\end{eqnarray}
Let $v_1(t)$ and $v_2(t)$ be the solutions of
\begin{eqnarray}
\label{mod20}
\dot{v}_1(t)&=& -\theta(d+\delta)\delta_2
-\gamma v_1(t), \quad t\neq n\tau,\nonumber\\
v_1(t^+)&=& v_1(t)+v_i, \quad t=n\tau,\nonumber
\end{eqnarray}
\mbox{and}
\begin{eqnarray}
\label{mod21}
\dot{v}_2(t)&=& \theta(d+\delta)\delta_2
-\gamma v_2(t), \quad t\neq n\tau,\nonumber\\
v_2(t^+)&=& v_2(t)+v_i, \quad t=n\tau,\nonumber
\end{eqnarray}
respectively. Then, the solution will be
\begin{eqnarray}\label{mod22}
v_1^*(t) &=& \frac{v_i e^{-\gamma (t-n\tau)}}{1
-e^{-\gamma \tau}}+\theta(d+\delta)\delta_2,\nonumber\\
v_2^*(t) &=& \frac{v_i e^{-\gamma (t-n\tau)}}{1
-e^{-\gamma \tau}}-\theta(d+\delta)\delta_2.
\end{eqnarray}

From (\ref{mod22}), it is clear that when $\delta_2\rightarrow 0$ we have
$v_1^*(t)\rightarrow v^*(t)$ and $v_2^*(t)\rightarrow v^*(t)$. Hence, it follows
from (\ref{mod19}) that $v(t)\rightarrow v^*(t)$ as $t\rightarrow \infty$.

Similarly, we can choose $0<\delta_3<\mu $ and, in the same way,
we can prove that $s(t)\rightarrow s^*(t)$ as $t\rightarrow \infty$.

Finally, we shall prove that $x(t)\rightarrow k$ as $t\rightarrow \infty$.
We already proved that $y(t),~z(t) \rightarrow 0$ as $t\rightarrow \infty$.
Thus, for $\delta_3 > 0$, there exists $t_3>0$ such that $y(t),~z(t) < \delta_3$
for all $t>t_3$. Hence, from the first equation of system (\ref{mod1}), we can write that
\begin{eqnarray}
\label{mod23}
r x(t)-\frac{r x^2(t)}{k}>\dot{x}(t)
\geq \{r-\delta_3\alpha(1+\phi)\} x(t)-\frac{r x^2(t)}{k},\nonumber
\end{eqnarray}
which implies that
\begin{eqnarray}
\label{mod24}
\frac{k x_0}{x_0+(k-x_0)e^{r t}} \geq x(t)
\geq\frac{k \{r-\delta_3\alpha(1+\phi)\}x_0}{r x_0
+[k\{r -\delta_3\alpha(1+\phi)\}
-r x_0]e^{-\{r-\delta_3\alpha(1+\phi)\} t}}.
\end{eqnarray}
Hence, for $\delta_3 \rightarrow 0$, $x(t)\rightarrow k$ as $t\rightarrow\infty$.
Thus, for application of biopesticide and chemical pesticide together
with the same time interval $t=n\tau$, we can say that system (\ref{mod1})
is locally as well as globally stable if
\begin{eqnarray}
\label{cond2}
&& c_1 \alpha -d-\frac{\lambda v_i e^{-\gamma (t-n\tau)}}{1-e^{-\gamma \tau}}
- \frac{m_1 s_i e^{-\mu (t-n\tau)}}{1-e^{-\mu \tau}}<0, \nonumber\\
&& c_2 \phi k \alpha -(d+\delta) - \frac{m_2 s_i e^{-\mu (t-n\tau)}}{1-e^{-\mu \tau}}<0.
\end{eqnarray}

Two subcases arise here, namely,

Subcase I. Application of biopesticide with time interval $t{=n \tau_1}$.

In this case, $s_i=0$. Hence, system (\ref{mod1}) is locally as well as globally
stable around the periodic solution if
\begin{eqnarray}
\label{cond3b}
&& c_1 \alpha k-d-\frac{\lambda v_i
e^{-\gamma (t-n\tau_1)}}{1-e^{-\gamma \tau_1}} <0, \nonumber\\
&& c_2 \phi k \alpha -(d+\delta) <0.
\end{eqnarray}

Subcase II. Application of chemical
pesticide with time interval $t{=n \tau_2}$.

In this case, $v_i=0$, and hence system (\ref{mod1}) is locally
as well as globally stable around the periodic solution if
\begin{eqnarray}
\label{cond4}
&& c_1 \alpha k-d- \frac{m_1 s_i
e^{-\mu (t-n\tau_2)}}{1-e^{-\mu \tau_2}}<0, \nonumber\\
&& c_2 \phi k \alpha -(d+\delta)
- \frac{m_2 s_i e^{-\mu (t-n\tau_2)}}{1-e^{-\mu \tau_2}}<0.
\end{eqnarray}
The proof is complete.
\end{proof}


\section{Numerical Simulations}
\label{sec4}

Now we solve the impulsive system numerically and we graphically display the results found.
We varied the crucial parameters within their feasible ranges to observe their impact
on the impulsive model's solution trajectories and equilibria. Precisely, we solved the
impulsive system and plotted the results in figures using the \texttt{ode45} MATLAB solver.

\begin{figure}[H]
\includegraphics[scale=0.7]{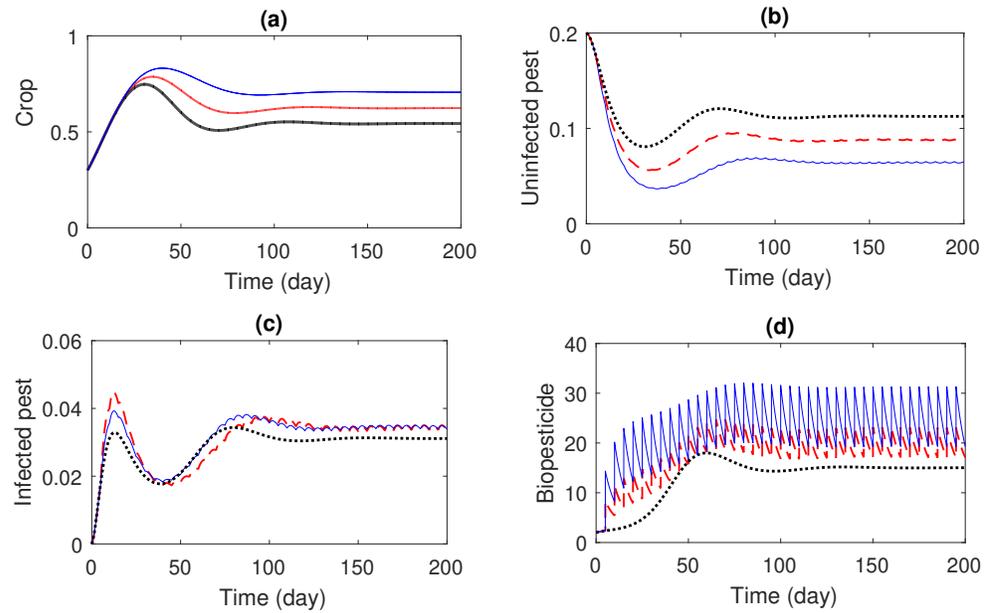}
\caption{Impact of biopesticide application in impulsive mode on system \eqref{mod1}.
Evolution of (\textbf{a}) crop; (\textbf{b}) uninfected pest; 
(\textbf{c}) infected pest; (\textbf{d}) biopesticide.
The set of parameters are $r = 0.1$, $k=1$, $\alpha=0.2$, $\beta=0.003$,
$m_1=0.8$, $m_2=0.6$, $c_1=0.5$, $c_2=0.8$, $\gamma=0.15$, $\delta=0.2$,
$d=0.05$, $\kappa=100$, $s=0.3$, and $\lambda=0.35$. Here, the time interval
is $\tau_1=5$ days and the rates of biopesticide release are $v_i=0$ (black line),
$v_i=6$ (red line), and $v_i=12$ (blue line).}\label{Fig1}
\end{figure}

In Figure~\ref{Fig1}, the impulsive time interval for microbial biological
pesticide release is 5~days, and releasing of biopesticide was considered
in different rates: $v_i=0$ (i.e., without pesticides), $v_i=6$, and $v_i=12$.
It is revealed that susceptible pest population decreases with an increase
in the release rate of biopesticides.

\begin{figure}[H]
\includegraphics[scale=0.7]{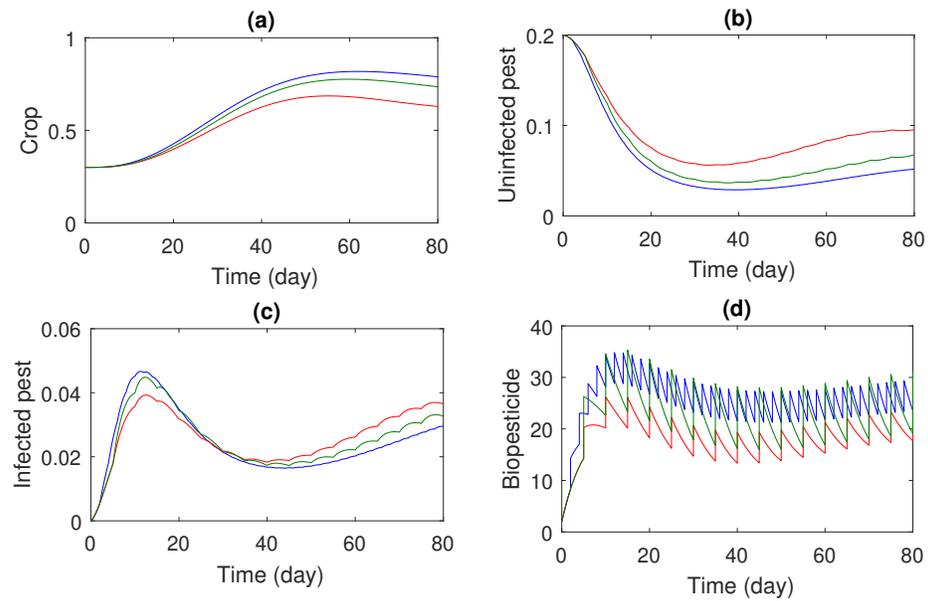}
\caption{Impact of biopesticide on system \eqref{mod1} for different impulse intervals
and rates.  Evolution of (\textbf{a}) crop; 
(\textbf{b}) uninfected pest; (\textbf{c}) infected pest; 
(\textbf{d}) biopesticide.
Red line indicates $v_i=6$ and $\tau_1=5$,
green line indicates $v_i=12$ and $\tau_1=5$,
and blue line indicates $v_i=12$ and $\tau_1=2$.}\label{Fig2}
\end{figure}

In Figure~\ref{Fig2}, by taking different impulsive intervals, biopesticide
is applied to the system. A better result is obtained for lower intervals (2 days)
but, with a higher release of biopesticide, pests are present in the system.
From Figures~\ref{Fig1} and \ref{Fig2}, we can conclude that pest control
using only biopesticides is very costly and a time-consuming process.

\begin{figure}[H]
\includegraphics[scale=0.7]{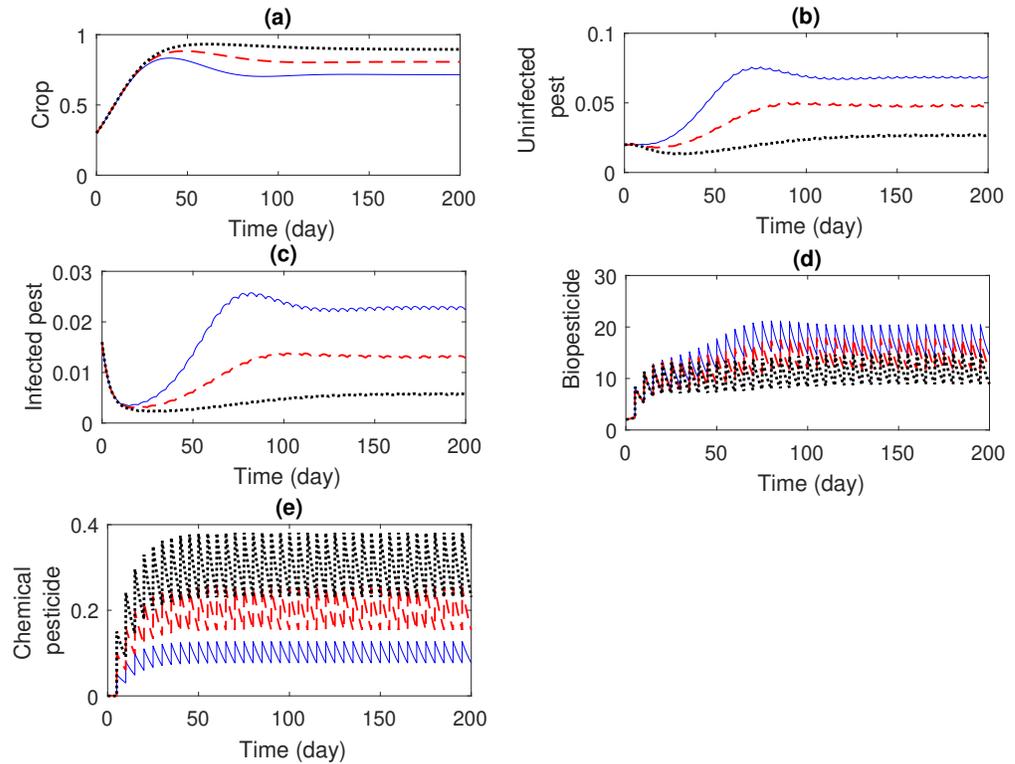}
\caption{Impact of both biopesticide and chemical pesticide
on system \eqref{mod1} with the same impulse interval, $\tau_1=\tau_2=5$.
Evolution of (\textbf{a}) crop; (\textbf{b}) uninfected pest; 
(\textbf{c}) infected pest; (\textbf{d}) biopesticide; (\textbf{e})~chemical pesticide.
The rates of impulses are $s_i=0.15$ and $v_i=6$ for black dotted color;
$s_i=0.1$ and $v_i=6$ for red dashed line;
and $s_i=0.05$ and $v_i=6$ for blue solid line.}\label{Fig3}
\end{figure}

\begin{figure}[H]
\includegraphics[scale=0.7]{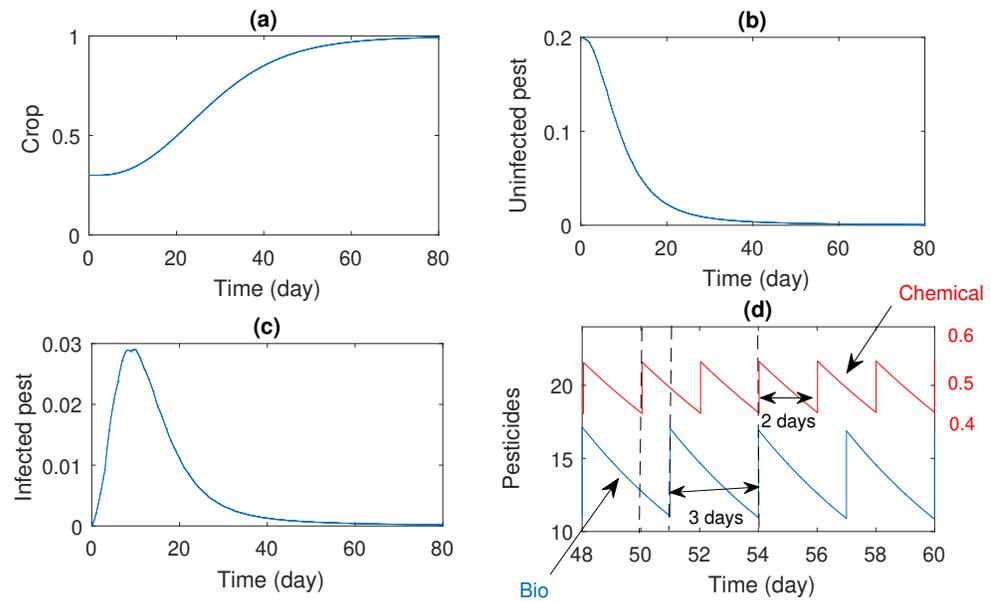}
\caption{Impact of both biopesticide and chemical pesticide
on system \eqref{mod1} with same impulse interval, $\tau_1=3$, $\tau_2=2$,
and where the rates of impulses are $s_i=0.15$ and $v_i=6$.
Evolution of (\textbf{a})~crop; (\textbf{b}) uninfected pest; (\textbf{c}) infected pest;
(\textbf{d}) pesticides (biopesticide in blue, chemical pesticide in red).}\label{Fig4}
\end{figure}

Recall that in our model we take $\tau_1$ as the time period for
biopesticide spraying (generally a virus particle) and $\tau_2$
as the time period for chemical pesticide sprays. In Figure~\ref{Fig3}
we see the effect for the same time intervals, $\tau_1=\tau_2=5$ days,
whereas in Figure~\ref{Fig4} we see the effect for different
time intervals, $\tau_1=3$ days and $\tau_2=2$ days.

If both microbial biopesticides and chemical pesticides are released simultaneously,
with an equal time interval of 5 days, then the extinction of both infected
and susceptible pest populations is possible (see Figure~\ref{Fig3}).

Figure~\ref{Fig4} illustrates the dynamics of the double impulse with different
impulse intervals. Double impulses occur at the time which is the common multiple
of the two intervals. For example, if we take $\tau_1=2$ and $\tau_2=3$,
then simultaneous impulses will occur at the times $t=6$, $t = 12$, $t=18$, and so on.
Figure~\ref{Fig4} is the most important figure characterizing the impact of two different
but simultaneous impulses on the total pest population with different time intervals.
It is observed that for $v_i=12$, $s_i=0.15$, $\tau_1=3$, and $\tau_2=2$, the total
pest population becomes extinct. In Figure~\ref{Fig4}d, the impact of double impulses
occurs at $t=48$, $t = 54$, $t=60$ days, etc., which are common multiples
of the two intervals $\tau_1=2$ and $\tau_2=3$.

It is also numerically checked that when the rate of the impulse control is high,
a comparatively lower interval can be taken for cost-effectiveness of the process.

Thus, the advantage of the impulsive control is that we can determine
the proper rate and a suitable interval of giving controlling agents
in the system.


\section{Discussion and Conclusions}
\label{sec5}

In the present research, we studied impulsive periodic applications
of integrated pesticides, that is, simultaneous use of biopesticide
and chemical pesticide in a pest management system. We proposed
a two-impulse mathematical model using an impulsive differential equation
to observe the impact of periodic application of the combined pesticides
in impulsive modes.

In the previous models available in the literature,
chemical and biological pesticides were used
in the model in a continuous way. In contrast, here, we used them
in an impulsive periodic way. Consequently, a two-impulse mathematical model
was established. Moreover, in the proposed model, we took
chemical pesticides concentration as the model population,
which is a novel approach.

Stability theory (Floquet theory) and numerical calculations were used
to examine the system dynamical behavior. We determined the conditions
under which the impulsive system will be stable both locally and globally. For example,
the local stability of a pest-free periodic orbit was established. The dynamics varied
with the rate of both biopesticide recruitment and the chemical pesticide concentration.

Chemical pesticides minimize the oscillations in the system and make the system stable
in a shorter time. Numerical and analytical analysis reveals that increasing frequency
of pesticide application will require less administration of biopesticide and chemical
pesticide, which is economically beneficial and environmentally safe.

Our research is directed to optimize and find the right combination of pesticides
with maximum benefit to the crop plant. The numerical simulation also shows that control
over the spraying of chemical pesticides is needed to control pests and minimize the cost
of cultivation. On the other hand, chemical pesticides may have negative environmental
implications due to their lingering effects; nonetheless, the best control approach
provides the least amount of collateral damage to the environment.

In a nutshell, the promising feature of the system is the combined use of the pesticides
in impulsive control methods that reduce the cost and negative effects on the environment.
Using a combination of pesticides to deliver the pesticide can save the cost and reduce
the side effects of chemical pesticides. Our obtained results will give a new perspective to
farmers who implement this in a real-world setting.

In the future, one can extend this work to an optimal impulsive system
for cost-effectiveness of the control process.


\vspace{6pt}

\authorcontributions{Conceptualization, F.A.B. and J.C.;
methodology, F.A.B., J.C. and D.F.M.T.;
software, F.A.B. and J.C.;
validation, F.A.B., J.C. and D.F.M.T.;
formal analysis, F.A.B., J.C. and D.F.M.T.;
investigation, F.A.B., J.C. and D.F.M.T.;
writing---original draft preparation, F.A.B., J.C. and D.F.M.T.;
writing---review and editing, F.A.B., J.C. and D.F.M.T.;
visualization, F.A.B. and J.C.
All authors have read and agreed
to the published version of the manuscript.}

\funding{Torres was funded by
The Portuguese Foundation for Science and Technology (FCT)
grant number UIDB/04106/2020.}

\institutionalreview{Not applicable.}

\informedconsent{Not applicable.}

\dataavailability{No new data were created or analyzed in this study.
Data sharing is not applicable to this article.}

\acknowledgments{The authors are very grateful
to four anonymous reviewers for careful reading
of the submitted manuscript and also for providing several
important comments and suggestions.}

\conflictsofinterest{The authors declare no conflicts of interest.
The funders had no role in the design of the study; in the collection,
analyses, or interpretation of data; in the writing of the manuscript;
or in the decision to publish the~results.}


\begin{adjustwidth}{-\extralength}{0cm}

\reftitle{References}



\PublishersNote{}

\end{adjustwidth}

\end{document}